\documentclass[11pt]{amsart}

\usepackage{amsmath}
\usepackage{amsthm}
\usepackage{amssymb}
\usepackage{amscd}

\newcommand{\gothic}{\mathfrak}
\newcommand{\p}{{\gothic{p}}}
\newcommand{\q}{{\gothic{q}}}

\newcommand{\m}{{\gothic{m}}}
\newcommand{\n}{{\gothic{n}}}

\newcommand{\Spec}{\operatorname{Spec}}

\newcommand{\Ext}{\operatorname{Ext{}}}

\newcommand{\Ann}{\operatorname{Ann}}

\renewcommand{\hat}{\widehat}

\renewcommand{\phi}{\varphi}

\renewcommand{\to}{{\longrightarrow}}

\newtheorem{thm}{Theorem}

\newtheorem{prop}[thm]{Proposition}
\newtheorem{lemma}[thm]{Lemma}
\newtheorem{defn}[thm]{Definition}
\newtheorem{conj}[thm]{Conjecture}

\begin{document}

\title{Local Rings of Countable Cohen-Macaulay Type }

\author{Craig Huneke}

\address{Department of Mathematics \\ 
        University of Kansas \\ 
        Lawrence, KS
        66045} 

\email{huneke@math.ukans.edu}

\urladdr{http://www.math.ukans.edu/\textasciitilde huneke}

\author{Graham J. Leuschke}


\email{gleuschke@math.ukans.edu}

\urladdr{http://www.math.ukans.edu/\textasciitilde gleuschke}

\date{\today}

\thanks{Both authors were supported by the National Science Foundation.}

\bibliographystyle{amsplain}

 \numberwithin{thm}{section}

\begin{abstract} We prove (the excellent case of) Schreyer's
conjecture that a local ring with countable CM type has at most a
one-dimensional
singular locus. Furthermore we prove that the localization of a
Cohen-Macaulay
local ring of countable CM type is again of countable CM type.
\end{abstract}  

\maketitle

Let $(R, \m)$ be a (commutative Noetherian) local ring of dimension
$d$.  Recall that a nonzero $R$-module $M$ is called {\it maximal
Cohen--Macaulay} (MCM) provided it is finitely generated and there
exists an $M$-regular sequence $\{x_1, \ldots, x_d\}$ in the maximal
ideal $\m$.  We say that $R$ itself is Cohen--Macaulay (CM) if it is MCM
as a module over itself.

The CM local rings of finite CM-representation type (meaning that they
have only finitely many nonisomorphic indecomposable MCM modules) have
been carefully studied over the last twenty years.  The complete
equicharacteristic hypersurfaces of finite CM type have been completely
classified (\cite{GreuelKnorrer}, \cite{BGS}, \cite{Knorrer}), as have
the complete equicharacteristic 2-dimensional normal domains
(\cite{Auslander:rationalsing}).  More generally, it is known that a CM
local ring of finite CM type has at most an isolated singularity (proved
by Auslander \cite{Auslander:isolsing} in the complete case,
Leuschke-Wiegand \cite{Leuschke-Wiegand:2000} in the excellent case, and
Huneke-Leuschke \cite{Huneke-Leuschke:2002} in general).  Yoshino's
monograph \cite{Yoshino:book} is a comprehensive source for
information about rings of finite CM type.

The related property of {\it countable} CM type has received much less
attention.  Buchweitz, Greuel, and Schreyer \cite{BGS} classified the
complete hypersurface singularities of countable CM type, but very
little more has been learned since then.  

The open questions and conjectures in Schreyer's 1987 survey article
\cite{Schreyer:1987} have inspired work on both finite and
countable CM-representation type.  For example, Conjecture~7.3(a) states
that a CM local ring $R$ has finite CM type if and only if the $\m$-adic
completion has finite CM type; this was recently proved in case $R$ is
excellent in \cite{Leuschke-Wiegand:2000}.  This paper is concerned with
another of Schreyer's conjectures:

\begin{conj}[\cite{Schreyer:1987}]\label{Schreyer} An analytic local ring
over the complex numbers of
countable CM type has at most a one-dimensional singular locus, that is,
$R_\p$ is regular for all primes $\p$ with $\dim R/\p > 1$.\end{conj}

We verify Conjecture~\ref{Schreyer} more generally for all excellent CM local rings
satisfying countable prime avoidance (Lemma~\ref{avoid}).
  Some assumption of uncountability is necessary
to avoid the degenerate case
of a countable ring, which {\it a fortiori}  has
only countably many isomorphism classes of modules.

\bigskip
\section{ Schreyer's Conjecture}
\medskip

\begin{defn} A Cohen-Macaulay local ring $(R, \m)$ is said to have {\rm
finite ({\it resp.},
countable) Cohen--Macaulay type} if it has only finitely (resp.,
countably) many isomorphism classes of maximal Cohen--Macaulay modules.
\end{defn}

For the proof of Schreyer's conjecture, we need the following well-known
lemma.

\begin{lemma}[countable prime avoidance
{\cite[Lemma 3]{Burch:1972}}; see also \cite{Sharp-Vamos}]\label{avoid}
 Let $A$ be a Noetherian local ring
which either is complete or has uncountable residue field.  Let
$\{\p_i\}$, $i=1,2,\ldots$, be a countable family of prime ideals of
$A$,
${\mathfrak a}$ an ideal of $A$, $x\in A$. Then $x+\mathfrak a\subseteq
\p_j$ for some $j$ whenever $x+\mathfrak
a\subseteq\bigcup^\infty_{i=1}\p_i$.\end{lemma}

\begin{thm}\label{mainthm}  Let $(R,\m)$ be an excellent Cohen-Macaulay
local ring of
dimension $d$, and assume either that $R$ is complete or that the
residue field $R/\m$ is uncountable. If $R$ has countable CM type, then
the singular locus of $R$ has dimension at most one.\end{thm}

\begin{proof} Assume that the singular locus of $R$ has dimension
greater than one.  Since $R$ is excellent, this means that the singular
locus is defined by an ideal $J$ of height strictly less than $d-1$. 
Let
$\{M_i\}_{i=1}^\infty$ be a complete list of representatives for the
isomorphism classes of indecomposable MCM $R$-modules.  Consider the set 
$$\Lambda = \{\p \in \Spec(R) \ |\ \p = \Ann(\Ext_R^1(M_i,M_j)),
\text{for some }i, j,\text{ and }\dim(R/\p)=1\}.$$
Note that $\Lambda$ is at most countable, and that $J$ is contained in
each $\p
\in \Lambda$.  By countable prime avoidance (applied to $R/J$), the
maximal ideal $\m$ is not contained in the union of all $\p$ in
$\Lambda$, so there is an element $f \in \m \setminus \bigcup_{\p \in
\Lambda}\p$.  Choose a prime $\q$ containing $f$ and $J$ such that $\dim
R/\q=1$; then of course $\q \not\in \Lambda$. 

Let $X$ (resp. $Y$) be a $(d-1)^{\text{th}}$ (resp. $d^{\text{th}}$)
syzygy
of $R/\q$.  Then $X$ and $Y$ are both MCM $R$-modules and we have a
nonsplit short exact sequence
\begin{equation}\label{syzygy}\tag{*} 0 \to Y \to F \to X \to 0, 
\end{equation}
where $F$ is a free $R$-module.  We claim that
$\Ann(\Ext_R^1(X,Y))=\q$.  It is clear that $\q$ kills $\Ext_R^1(X,Y)
\cong \Ext_r^{d+1}(R/\q,Y)$.  To see the opposite containment, note that
since $\q$ contains $J$, $R_\q$ is not regular.  The resolution of the
residue field of $R/\q$ is thus infinite, and neither $X_\q$ nor $Y_\q$
is free, so (\ref{syzygy}) is nonsplit when localized at $\q$.

We can write both $X$ and $Y$ as direct sums of copies of the
indecomposables $M_i$, and further write
\begin{equation*}\Ext_R^1(X,Y) \cong \bigoplus_{i,j}
\Ext_R^1(M_i,M_j)^{a_{ij}}\end{equation*}
with all but finitely many of the $a_{ij}$ equal to zero.  Then $\q$ is
the intersection of the annihilators of the nonzero Ext modules
appearing
in the above decomposition. Since $\q$ is prime, it must equal one of
these
annihilators, and then $\q \in
\Lambda$, a contradiction.\end{proof}

\bigskip
\section{Localization of rings with countable CM type}
\medskip

Let $(R,\m)$ be an excellent local ring of countable CM type, and assume
$R$ either
has an uncountable residue field or is complete. By
Theorem~\ref{mainthm},
the dimension of the singular locus of $R$ is at most one. Thus there  
are at most finitely many prime ideals $\p_1,\ldots,\p_n$ such that
$R_{\p_i}$ is
not regular and $\p_i\ne \m$. All such primes have dimension one, i.e.,
dim$(R/\p_i) = 1$ for $i = 1,\ldots,n$. To understand the structure of
these
rings, one wishes to know what type of singularity $R$ has at these
primes.
A quick inspection of the list of examples given in 
\cite{Schreyer:1987} shows that $R_{\p_i}$ has finite
CM type! Our next main result shows
that countable CM type localizes; hence in general each $R_{\p_i}$ has
countable
CM type. Note also that in the case $R$ is complete, while $R_{\p_i}$
is no longer complete in general, it will have uncountable residue
field.

\begin{thm} Let $(R,\m)$ be a CM local ring with a canonical module.  If
$R$ has countable CM type, then $R_\p$ has countable CM type for each
prime $\p$ of $R$.\end{thm}

\begin{proof} Let $\p \in \Spec(R)$ and assume that $\{M_\alpha\}$ is an
uncountable family of $R$-modules such that $\{(M_\alpha)_\p\}$ are
nonisomorphic indecomposable MCM $R_\p$ modules.  For each $\alpha$
there is a {\it Cohen--Macaulay approximation} of $M_\alpha$, that is, a
short exact sequence
\begin{equation}\tag*{($\chi_\alpha$)}\label{CMapprox}0 \to Y_\alpha \to
X_\alpha \to M_\alpha \to 0\end{equation}
wherein $X_\alpha$ is a MCM $R$-module and $Y_\alpha$ has finite
injective dimension \cite{Auslander-Buchweitz}.

Since there are uncountably many modules $X_\alpha$, there must be
uncountably many $X_\alpha$ of some fixed multiplicity. Fixing that
multiplicity,  and using that there are only countably many
isomorphism classes of MCMs, we then find that there are uncountably
many short exact sequences
\begin{equation}\tag*{($\chi_\beta$)}\label{CMapproxx}0 \to Y_\beta \to
X \to M_\beta \to 0\end{equation}
where $X$ is a fixed MCM $R$-module, $Y_\beta$ has finite injective
dimension,
and the $M_\beta$ are among our original list of $M_\alpha$.

Since each $(M_\beta)_\p$ is a MCM $R_\p$-module and $(Y_\beta)_\p$
has finite injective dimension over $R_\p$,
$\Ext_R^1(M_\beta,Y_\beta)_\p \cong \Ext_R^1((M_\beta)_\p,
(Y_\beta)_\p) = 0$. This follows from
\cite[Proposition 4.9]{Ro}: if $Y$ is a finitely
generated module having finite injective dimension, then for all
finitely generated $R$-modules $M$,
$$\text{depth}(M)  + \text{sup}\{i|\, \text{Ext}^i_R(M,Y)\ne 0\} =
\text{depth}(R).$$

In particular, each
extension \ref{CMapproxx} splits when localized at $\p$.  This implies
that $(M_\beta)_\p$ is a direct summand of $(X)_\p$ for each
$\beta$. But over a local ring, a finitely generated  $R$-module $Q$
can have at most finitely many non-isomorphic summands.
\footnote{This can be seen by passing to the completion, where the
Krull-Schmidt
theorem holds.  After completion, any direct summand of $\hat{Q}$ must
be isomorphic to a direct sum of a fixed finite subset of  the
indecomposable summands of the completion of $Q$.  Hence there are
only finitely many  such isomorphism classes after completion,
and since the completions of two modules  are isomorphic if and only
if the  two modules are isomorphic, we are done (see \cite[Proposition
(2.5.8) Chap. IV, Section 2]{Gr}).}
Since there are uncountably many $(M_\beta)_\p$ which must
be summands of $X_\p$, this contradiction proves the theorem.
\end{proof}

The results above, together with known examples, suggest a
plausible question:
\smallskip

{\it Let $R$ be a complete local Cohen-Macaulay ring of countable
CM type, and assume that $R$ has an isolated singularity. Is $R$
then necessarily of finite CM type?}

\smallskip
We end the paper with an observation that having countable CM type
descends from faithfully flat overrings. The method follows that of
\cite{Wiegand:1998}.
Countable CM type cannot in general ascend to the completion, since
countable rings are of countable CM type, but their completions
are not of countable CM type if the dimension of the singular locus
is at least two. However, when the residue field is uncountable,  we
do not know if countable CM type
ascends to the completion.
The proof of the analogous assertion for finite CM type 
uses that the ring is necessarily Gorenstein on the punctured spectrum,
which we do not know is true for countable CM type.

\begin{prop}Let $(R, \m)$ be a CM local ring and $(S, \n)$ a faithfully
flat $R$-algebra such that the closed fibre $S/\m S$ is CM.  If $S$ has
countable CM type, then so does $R$.\end{prop}

\begin{proof}  Since the closed fibre is CM, the extension $S\otimes_R
M$ of any MCM $R$-module is an MCM $S$-module.  Let $Z_1, Z_2, \ldots$
be a complete list of all indecomposable MCM $S$-modules which are
direct summands of $S\otimes_R M$ for some MCM $R$-module $M$ (ignore
any $S$-modules that do not appear in a direct-sum decomposition of some
extended module).  For each $i=1, 2, \ldots$, choose an indecomposable
MCM $R$-module $X_i$ so that $Z_i \oplus W_i \cong S \otimes_R X_i$ for
some $S$-module $W_i$.

For an indecomposable MCM $R$-module $N$, write $S \otimes_R N \cong
\bigoplus_i Z_i^{a_i}$, where all but finitely many of the $a_i$ are
zero.  We assume that $a_i =0$ for $i>n$ and write the sum as a finite
one.  Then
$$(S\otimes_R N) \oplus W_1^{a_1} \oplus \cdots \oplus W_n^{a_n} \cong
S\otimes_R(X_1^{a_1}\oplus \cdots \oplus X_n^{a_n}),$$
so $S\otimes_R N$ is a direct summand of $S\otimes_R(X_1^{b}\oplus
\cdots \oplus X_n^{b})$, where $b = \max\{a_i\}$.  In other words,
$S\otimes_R N$ is in the ``plus category'' of $S\otimes_R(X_1\oplus
\cdots \oplus X_n)$ (see \cite{Wiegand:1998}).  By \cite[Lemma
1.2]{Wiegand:1998}, $N$ is in the plus category of $X_1\oplus \cdots
\oplus X_n$, and by \cite[Theorem 1.1]{Wiegand:1998}, there are only
finitely many possible such $N$.  Since the set of all finite subsets of
$\{X_1, X_2, \ldots\}$ is a countable set, this shows that $R$ has only
countably many indecomposable MCM modules up to isomorphism.\end{proof}

\providecommand{\bysame}{\leavevmode\hbox to3em{\hrulefill}\thinspace}

\end{document}